\documentclass[a4paper,twoside,10pt]{amsart}
\usepackage{amsmath,amsfonts,amssymb,amsthm}
\newtheorem{theorem}{Theorem}[section]

\usepackage{mathrsfs}
\usepackage{color,graphicx}
\usepackage[a4paper]{geometry}
\numberwithin{equation}{section}

\author[G. Nemes]{Gerg\H{o} Nemes}
\address{Central European University, Department of Mathematics and its Applications, H-1051 Budapest, N\'ador utca 9, Hungary}
\email{nemesgery@gmail.com}

\keywords{asymptotic expansions, Barnes $G$-function, error bounds, Stokes phenomenon.}
\subjclass[2010]{41A60, 30E15, 34M40}

\begin{document}

\title[Asymptotics of the Barnes $G$-function]{Error bounds and exponential improvement\\ for the asymptotic expansion of\\ the Barnes $G$-function}

\begin{abstract} In this paper we establish new integral representations for the remainder term of the known asymptotic expansion of the logarithm of the Barnes $G$-function. Using these representations, we obtain explicit and numerically computable error bounds for the asymptotic series, which are much simpler than the ones obtained earlier by other authors. We find that along the imaginary axis, suddenly infinitely many exponentially small terms appear in the asymptotic expansion of the Barnes $G$-function. Employing one of our representations for the remainder term, we derive an exponentially improved asymptotic expansion for the logarithm of the Barnes $G$-function, which shows that the appearance of these exponentially small terms is in fact smooth, thereby proving the Berry transition property of the asymptotic series of the $G$-function.
\end{abstract}
\maketitle

\section{Introduction and main results}

In a series of papers published between 1899--1904 \cite{Barnes1,Barnes2,Barnes3,Barnes4}, Barnes introduced and studied certain generalisations of the classical Gamma function $\Gamma\left(z\right)$, called the multiple Gamma functions $\Gamma_k\left(z\right)$. Among them, the case $G\left(z\right) = 1/\Gamma_2\left(z\right)$, being the so-called double Gamma function or the Barnes $G$-function, is of special interest. This function has several important applications in many different areas of mathematics, such as the theory of $p$-adic $L$-functions \cite{Cassou-Nogues} or the study of determinants of Laplacians \cite{Choi,Osgood,Quine,Sarnak,Vardi}.

It is known that the $G$-function is an entire function and has the Weierstrass canonical product
\[
G\left( {z + 1} \right) = \left( {2\pi } \right)^{z/2} e^{ - \left[ {z\left( {z + 1} \right) + \gamma z^2 } \right]/2} \prod\limits_{n = 1}^\infty  {\left[ {\left( {1 + \frac{z}{n}} \right)^n e^{ - z + z^2 /\left( {2n} \right)} } \right]} ,
\]
where $\gamma=0.57721566\ldots$ is the Euler--Mascheroni constant. It also satisfies the functional equation $G\left( {z + 1} \right) = \Gamma \left( z \right)G\left( z \right)$.

In January 2000, Richard Askey proposed, in the panel discussion of the San Diego symposium on asymptotics and applied analysis, a deeper study of this
function. In particular, the derivation of an asymptotic expansion for $G\left(z\right)$ is necessary for approximating this function for large values of $z$.

Motivated by the function's increasing interest and Askey's proposition, Ferreira and L\'{o}pez \cite{Ferreira1} established the following expansion in terms of the Bernoulli numbers
\begin{gather}\label{eq2}
\begin{split}
\log G\left( {z + 1} \right) = \frac{1}{4}z^2 & + z\log \Gamma \left( {z + 1} \right) - \left( {\frac{1}{2}z\left( {z + 1} \right) + \frac{1}{{12}}} \right)\log z - \log A \\ & + \sum\limits_{n  = 1}^{N - 1} {\frac{{B_{2n  + 2} }}{{2n \left( {2n  + 1} \right)\left( {2n  + 2} \right)z^{2n } }} } + R_N\left( z \right),
\end{split}
\end{gather}
for $\left| {\arg z} \right| < \frac{\pi}{2}$, where
\begin{equation}\label{eq1}
R_N\left(z\right) = \int_0^{+\infty} {\left( {\frac{t}{2}\coth \left( {\frac{t}{2}} \right) - \sum\limits_{k = 0}^N {\frac{{B_{2k} }}{{\left( {2k} \right)!}}t^{2k} } } \right)\frac{{e^{ - zt} }}{{t^3 }}dt}.
\end{equation}
The constant $A$ is called the Glaisher--Kinkelin constant, and can be expressed as
\[
\log A = \frac{{\gamma  + \log \left( {2\pi } \right)}}{{12}} - \frac{{\zeta '\left( 2 \right)}}{{2\pi ^2 }} = \frac{1}{{12}} - \zeta '\left( { - 1} \right) = 0.24875447\ldots,
\]
where $\zeta$ denotes Riemann's Zeta function. They also showed that this expansion can be extended to the wider sector $\left| {\arg z} \right| < \pi$, by rotating the path of integration in \eqref{eq1}. By estimating the error term $R_N\left(z\right)$, they deduced the following asymptotic expansion for the Barnes $G$-function:
\begin{gather}\label{eq4}
\begin{split}
\log G\left( {z + 1} \right) \sim \frac{1}{4}z^2  + z\log \Gamma \left( {z + 1} \right) - \left( {\frac{1}{2}z\left( {z + 1} \right) + \frac{1}{{12}}} \right)\log z - \log A \\ + \sum\limits_{n  = 1}^{\infty} {\frac{{B_{2n  + 2} }}{{2n \left( {2n  + 1} \right)\left( {2n  + 2} \right)z^{2n } }}} ,
\end{split}
\end{gather}
as $z\to \infty$ in the sector $\left| {\arg z} \right| \leq \pi -\delta< \pi$ with any fixed $0<\delta\leq \pi$.

In this paper we use \eqref{eq2} as the definition of $R_N\left( z \right)$ for $\left| {\arg z} \right| < \pi$ and $N\geq 1$. Throughout the paper, empty sums are taken to be zero.

The aim of this paper is to establish new integral representations for the remainder $R_N\left( z \right)$. Using these representations, we obtain several new properties of the asymptotic expansion \eqref{eq4}, including simple, explicit and numerically computable bounds for $R_N\left( z \right)$, exponentially improved asymptotic expansion, and the smooth transition of the Stokes discontinuities.

In our first theorem, we give two integral representations for the remainder $R_N\left( z \right)$ of the asymptotic series \eqref{eq4}. In formula \eqref{eq31} below, $\mathop{\text{Li}_2} \left( w \right)$ denotes the Dilogarithm, defined by
\begin{equation}\label{eq18}
\mathop{\text{Li}_2} \left( w \right) = \sum\limits_{n = 1}^\infty  {\frac{{w^n }}{{n^2 }}} 
\end{equation}
for $\left|w\right| \leq 1$. 

\begin{theorem}\label{thm1} For any $N\geq 1$, the remainder $R_N \left( z \right)$ in \eqref{eq2} has the representations
\begin{align}
R_N \left( z \right) & = \frac{1}{{z^{2N} }}\frac{{\left( { - 1} \right)^{N + 1} }}{\pi }\int_0^{ + \infty } {\left( {\int_0^1 {\frac{{s^{2N - 1} }}{{1 + \left( {st/z} \right)^2 }}ds} } \right)t^{2N} \log \left( {1 - e^{ - 2\pi t} } \right)dt} \label{eq5} \\
& = \frac{1}{{z^{2N} }}\frac{{\left( { - 1} \right)^N }}{{2\pi ^2 }}\int_0^{ + \infty } {\frac{{t^{2N - 1} }}{{1 + \left( {t/z} \right)^2 }}\mathop{\mathrm{Li}_2} \left( {e^{ - 2\pi t} } \right)dt} \label{eq31}
\end{align}
for $\left|\arg z\right|<\frac{\pi}{2}$, and
\begin{align}
R_N \left( z \right) & =  - \frac{1}{{2N\left( {2N + 1} \right)}}\int_0^{ + \infty } {\frac{{B_{2N + 1} \left( {t - \left\lfloor t \right\rfloor } \right)}}{{\left( {t + z} \right)^{2N} }}dt} \label{eq6} \\ &  =  - \frac{1}{{\left( {2N + 1} \right)\left( {2N + 2} \right)}}\int_0^{ + \infty } {\frac{{B_{2N + 2}  + B_{2N + 2} \left( {t - \left\lfloor t \right\rfloor } \right)}}{{\left( {t + z} \right)^{2N + 1} }}dt} \label{eq7}
\end{align}
for $\left|\arg z\right|<\pi$, where $B_n\left(x\right)$ denotes the $n$th Bernoulli polynomial.
\end{theorem}

We note that \eqref{eq31} can also be written in terms of the Bose--Einstein integral $G_p\left(w\right)$ \cite[pp. 611--612]{NIST} since $\mathop{\text{Li}_2} \left( {e^{ - 2\pi t} } \right) = G_1 \left( { - 2\pi t} \right)$. The second theorem provides two bounds for the error term $R_N\left( z \right)$, one being valid in the closed sector $\left|\arg z\right| \leq \frac{\pi}{2}$ the other being valid in the larger one $\left|\arg z\right| < \pi$.

\begin{theorem}\label{thm2} Let $\theta = \arg z$. For any $N\geq 1$, the remainder $R_N \left( z \right)$ in \eqref{eq2} satisfies
\begin{equation}\label{eq15}
\left| {R_N \left( z \right)} \right|  \le \frac{{\left| {B_{2N + 2} } \right|}}{{2N\left( {2N + 1} \right)\left( {2N + 2} \right)\left| z \right|^{2N} }} \begin{cases} \min\left(\left|\csc\left(2\theta\right)\right|,\sqrt {\frac{e}{4}\left( {2N + \frac{5}{2}} \right)} \right) & \text{if } \; \frac{\pi}{4} < \left|\theta\right| \leq \frac{\pi}{2} \\ 1 & \text{if } \; \left|\theta\right| \leq \frac{\pi}{4} , \end{cases}
\end{equation}
and
\begin{equation}\label{eq17}
\left| {R_N \left( z \right)} \right| \le \frac{{\left| {B_{2N + 2} } \right|}}{{2N\left( {2N + 1} \right)\left( {2N + 2} \right)\left| z \right|^{2N} }}\sec ^{2N + 1} \left( {\frac{\theta }{2}} \right)
\end{equation}
for $\left|\arg z\right|<\pi$. Moreover, if $z>0$ then $R_N \left( z \right)$ has the sign of the first neglected term and its absolute value is smaller than that term.
\end{theorem}

In the third theorem, we give alternative bounds for the sectors $\frac{\pi}{4} < \left|\arg z\right| < \pi$.

\begin{theorem}\label{thm3} Let $\theta = \arg z$. For any $N\geq 1$, the remainder $R_N \left( z \right)$ in \eqref{eq2} satisfies
\[
\left| {R_N \left( z \right)} \right| \le \frac{{\csc \left( {2\left( {\theta  - \varphi^\ast } \right)} \right)}}{{\cos ^{2N + 1} \varphi^\ast }}\frac{{\left| {B_{2N + 2} } \right|}}{{2N\left( {2N + 1} \right)\left( {2N + 2} \right)\left| z \right|^{2N} }}
\]
for $\frac{\pi}{4}<\arg z <\pi$, where $0 < \varphi^\ast  < \frac{\pi}{2}$ is the unique solution of the implicit equation
\begin{equation}\label{eq16}
\left( {2N + 3 } \right)\cos \left( {3\varphi^\ast - 2\theta} \right) = \left( {2N - 1 } \right)\cos \left( {\varphi^\ast- 2\theta } \right),
\end{equation}
that satisfies $-\frac{\pi}{2}+\theta < \varphi^\ast <\frac{\pi}{2}$ if $\frac{3\pi}{4} \leq \theta <\pi$; $-\frac{\pi}{2}+\theta < \varphi^\ast <-\frac{\pi}{4}+\theta$ if $\frac{\pi}{2} \leq \theta <\frac{3\pi}{4}$; and $0 < \varphi^\ast <-\frac{\pi}{4}+\theta$ if $\frac{\pi}{4} < \theta <\frac{\pi}{2}$. Similarly, we have
\[
\left| {R_N \left( z \right)} \right| \le  - \frac{{\csc \left( {2\left( {\theta  - \varphi ^ *  } \right)} \right)}}{{\cos ^{2N + 1} \varphi ^ *  }}\frac{{\left| {B_{2N + 2} } \right|}}{{2N\left( {2N + 1} \right)\left( {2N + 2} \right)\left| z \right|^{2N} }}
\]
for $- \pi  < \arg z <  - \frac{\pi }{4}$, where $-\frac{\pi }{2} < \varphi^\ast  < 0$ is the unique solution of the equation \eqref{eq16}, that satisfies $ -\frac{\pi}{2} < \varphi^\ast < \frac{\pi}{2}+\theta$ if $-\pi <\theta \leq -\frac{3\pi}{4}$; $\frac{\pi}{4}+\theta < \varphi^\ast < \frac{\pi}{2}+\theta$ if $ -\frac{3\pi}{4} < \theta  \leq -\frac{\pi}{2}$; and $\frac{\pi}{4}+\theta < \varphi^\ast < 0$ if $-\frac{\pi}{2} < \theta <-\frac{\pi}{4}$.
\end{theorem}

Note that all these error bounds involve the absolute value of the first omitted term of the asymptotic expansion \eqref{eq4}, multiplied by a factor which depends only on $\arg z$ and the index $N$ of the first neglected term. The bounds obtained originally by Ferreira and L\'{o}pez \cite{Ferreira1} have more complicated forms. Later, employing a result of Pedersen \cite{Pedersen}, Ferreira \cite{Ferreira2} established the following simpler bound
\[
\left| {R_N \left( z \right)} \right| \le \frac{{\left| {B_{2N + 2} } \right|}}{{2N\left( {2N + 1} \right)\left( {2N + 2} \right)\left| z \right|^{2N} }}\sec ^{2N} \theta ,
\]
for $\left|\arg z\right| < \frac{\pi}{2}$ and $N\geq 1$. It is easy to see that this bound is weaker than \eqref{eq15}.

In our last theorem, we give a new convergent series expansion for the logarithm of the Barnes $G$-function. This expansion involves two infinite series. The terms of the first infinite series are finite sums whose length depend on an arbitrary sequence of non-negative integers $N_1,N_2,N_3,\ldots$. This infinite series is in fact a certain decomposition of the truncated version of the asymptotic series \eqref{eq4} (see the remark after the theorem). The terms of the second infinite series involve the Terminant function $\widehat T_p\left(w\right)$, whose definition and basic properties are given in Section \ref{section4}. In Subsection \ref{subsection41}, we will show that, as a manifestation of the Stokes phenomenon, suddenly infinitely many exponentially small terms appear in the asymptotic expansion \eqref{eq4} as $\arg z$ passes through the values $\pm \frac{\pi}{2}$. An optimal choice of the indices $N_1,N_2,N_3,\ldots$ and the asymptotic behaviour of the Terminant functions will show that the appearance of these exponentially small terms is in fact smooth, proving the Berry transition property of the asymptotic series \eqref{eq4} (see Subsection \ref{subsection42}). Moreover, with the optimal choice of the indices $N_1,N_2,N_3,\ldots$, \eqref{eq32} is an exponentially improved asymptotic expansion for $\log G\left( {z + 1} \right)$. We remark that a very similar behaviour was found for the asymptotic series of the Gamma function \cite{Berry3,Paris2}\cite[pp. 279--288]{Paris} and the Hurwitz Zeta function \cite{Paris3}.

\begin{theorem}\label{thm4} For an arbitrary sequence of non-negative integers $N_1,N_2,N_3,\ldots$, the logarithm of the Barnes $G$-function has the exact expansion
\begin{gather}\label{eq32}
\begin{split}
\log G\left( {z + 1} \right) = \frac{1}{4}z^2 & + z\log \Gamma \left( {z + 1} \right) - \left( {\frac{1}{2}z\left( {z + 1} \right) + \frac{1}{{12}}} \right)\log z - \log A \\ & - \sum\limits_{k = 1}^\infty  {\frac{1}{{\left( {2\pi k} \right)^2 }}\sum\limits_{n = 0}^{N_k  - 1} {\left( { - 1} \right)^n \frac{{2\left( {2n + 1} \right)!}}{{\left( {2\pi kz} \right)^{2n + 2} }}} } \\ & - \sum\limits_{k = 1}^\infty  {\left( {\frac{{\widehat T_{2N_k  + 1} \left( {2\pi kiz} \right)}}{{2\pi ik^2 }}e^{2\pi kiz}  + \frac{{\widehat T_{2N_k  + 1} \left( { - 2\pi kiz} \right)}}{{2\pi ik^2 }}e^{ - 2\pi kiz} } \right)},
\end{split}
\end{gather}
provided that $\left|\arg z\right| < \pi$.
\end{theorem}

The special case $N_k  = N - 1 \ge 0$ ($k\geq 1$) yields the following convergent series for the remainder $R_N\left( z \right)$ of the asymptotic series \eqref{eq4}:
\[
R_N \left( z \right) =  - \sum\limits_{k = 1}^\infty  {\left( {\frac{{\widehat T_{2N - 1} \left( {2\pi kiz} \right)}}{{2\pi ik^2 }}e^{2\pi kiz}  + \frac{{\widehat T_{2N - 1} \left( { - 2\pi kiz} \right)}}{{2\pi ik^2 }}e^{ - 2\pi kiz} } \right)} ,
\]
for $\left|\arg z\right| < \pi$.  In deriving this expansion, we used the well-known representation of the Bernoulli numbers in terms of the Riemann Zeta function \cite[equation 25.6.2]{NIST} and the definition of $R_N \left( z \right)$.

The rest of the paper is organised as follows. In Section \ref{section2}, we prove the formulas for the remainder term stated in Theorem \ref{thm1}. In Section \ref{section3}, we prove the error bounds given in Theorems \ref{thm2} and \ref{thm3}. In Section \ref{section4}, we prove the expansion presented in Theorem \ref{thm4}, and provide a detailed discussion of the Stokes phenomenon related to the asymptotic expansion \eqref{eq4}. The paper concludes with a short discussion in Section \ref{section5}.

\section{Proof of the representations for the remainder term}\label{section2}

In this section, we prove Theorem \ref{thm1}. In the first part we prove \eqref{eq5} and \eqref{eq31}, while in the second part we derive \eqref{eq6} and \eqref{eq7}.

\subsection{Proof of formulas \eqref{eq5} and \eqref{eq31}} Our starting point is the representation
\[
\log G\left( {z + 1} \right) = \frac{1}{2}z\log \left( {2\pi } \right) - \frac{1}{2}z\left( {z + 1} \right) + z\log \Gamma \left( {z + 1} \right) - \int_0^z {\log \Gamma \left( {s + 1} \right)ds} 
\]
for $\left|\arg z\right|<\pi$, with the logarithm taking its principal value \cite[equation 5.17.4]{NIST}. As for the path of integration, we take a straight line that connects $0$ to $z$. For the logarithm of the Gamma function $\log \Gamma \left( {s + 1} \right)$ under the integral, we introduce the expansion
\[
\log \Gamma \left( {s + 1} \right) = \left( {s + \frac{1}{2}} \right)\log  - s + \frac{1}{2}\log \left( {2\pi } \right) + \mu \left( s \right),
\]
where $\mu \left( s \right)$ is the Binet function and $\left|\arg s\right|<\pi$ \cite[pp. 56--58]{Remmert}. Integrating term-by-term, we deduce
\[
\log G\left( {z + 1} \right) = \frac{1}{4}z^2  + z\log \Gamma \left( {z + 1} \right) - \left( {\frac{1}{2}z\left( {z + 1} \right) + \frac{1}{{12}}} \right)\log z - \log A + R_1 \left( z \right)
\]
with
\[
R_1 \left( z \right) = \log A + \frac{1}{{12}}\log z - \int_0^z {\mu \left( s \right)ds} .
\]
Now suppose that $\left|\arg z\right|<\frac{\pi}{2}$ and $\left|\arg s\right|<\frac{\pi}{2}$. If we substitute the following formula
\[
\mu \left( s \right) =  - \frac{1}{\pi }\int_0^{ + \infty } {\frac{s}{{s^2  + t^2 }}\log \left( {1 - e^{ - 2\pi t} } \right)dt} ,
\]
valid for $\left|\arg s\right|<\frac{\pi}{2}$ \cite[p. 66]{Remmert}, we find
\begin{align*}
R_1 \left( z \right) & = \log A + \frac{1}{{12}}\log z + \frac{1}{{2\pi }}\int_0^{ + \infty } {\log \left( {1 + \left( {\frac{z}{t}} \right)^2 } \right)\log \left( {1 - e^{ - 2\pi t} } \right)dt} \\ & = \log A - \frac{1}{\pi }\int_0^{ + \infty } {\log t\log \left( {1 - e^{ - 2\pi t} } \right)dt}  + \frac{1}{{2\pi }}\int_0^{ + \infty } {\log \left( {1 + \left( {\frac{t}{z}} \right)^2 } \right)\log \left( {1 - e^{ - 2\pi t} } \right)dt} .
\end{align*}
The first integral can be calculated using integration by parts
\begin{align*}
 - \frac{1}{\pi }\int_0^{ + \infty } {\log t\log \left( {1 - e^{ - 2\pi t} } \right)dt} & = 2\int_0^{ + \infty } {\frac{{t\log t}}{{e^{2\pi t}  - 1}}dt}  - 2\int_0^{ + \infty } {\frac{t}{{e^{2\pi t}  - 1}}dt}  \\ & =  - \frac{{\gamma  + \log \left( {2\pi } \right)}}{{12}} + \frac{{\zeta '\left( 2 \right)}}{{2\pi ^2 }} =  - \log A,
\end{align*}
leading to the result
\[
R_1 \left( z \right) = \frac{1}{{2\pi }}\int_0^{ + \infty } {\log \left( {1 + \left( {\frac{t}{z}} \right)^2 } \right)\log \left( {1 - e^{ - 2\pi t} } \right)dt} ,
\]
which is valid for $\left|\arg z\right|<\frac{\pi}{2}$. Next, we substitute the expansion
\[
\log \left( {1 + \left( {\frac{t}{z}} \right)^2 } \right) = \sum\limits_{n = 1}^{N - 1} {\left( { - 1} \right)^{n + 1} t^{2n}\frac{1}{n z^{2n} }}  + \left( { - 1} \right)^{N + 1} t^{2N} \frac{1}{{z^{2N} }}\int_0^1 {\frac{{2s^{2N - 1} }}{{1 + \left( {st/z} \right)^2 }}ds} 
\]
with $N\geq 1$, to deduce that
\begin{align*}
R_1 \left( z \right) = \; & \sum\limits_{n = 1}^{N - 1} {\frac{1}{{2nz^{2n} }}\frac{{\left( { - 1} \right)^{n + 1} }}{\pi }\int_0^{ + \infty } {t^{2n} \log \left( {1 - e^{ - 2\pi t} } \right)dt} }  \\ & + \frac{1}{{z^{2N} }}\frac{{\left( { - 1} \right)^{N + 1} }}{\pi }\int_0^{ + \infty } {\left( {\int_0^1 {\frac{{s^{2N - 1} }}{{1 + \left( {st/z} \right)^2 }}ds} } \right)t^{2N} \log \left( {1 - e^{ - 2\pi t} } \right)dt} 
\\ = \; & \sum\limits_{n = 1}^{N - 1} {\frac{{B_{2n + 2} }}{{2n\left( {2n + 1} \right)\left( {2n + 2} \right)z^{2n} }}}  \\ & + \frac{1}{{z^{2N} }}\frac{{\left( { - 1} \right)^{N + 1} }}{\pi }\int_0^{ + \infty } {\left( {\int_0^1 {\frac{{s^{2N - 1} }}{{1 + \left( {st/z} \right)^2 }}ds} } \right)t^{2N} \log \left( {1 - e^{ - 2\pi t} } \right)dt} ,
\end{align*}
where we have used the known representation of the Bernoulli numbers \cite[equation 24.7.5]{NIST}
\begin{equation}\label{eq11}
\frac{{B_{2n + 2} }}{{\left( {2n + 1} \right)\left( {2n + 2} \right)}} = \frac{{\left( { - 1} \right)^{n + 1} }}{\pi }\int_0^{ + \infty } {t^{2n} \log \left( {1 - e^{ - 2\pi t} } \right)dt} .
\end{equation}
Finally, by the definition of $R_N \left( z \right)$, we conclude that
\[
R_N \left( z \right) = \frac{1}{{z^{2N} }}\frac{{\left( { - 1} \right)^{N + 1} }}{\pi }\int_0^{ + \infty } {\left( {\int_0^1 {\frac{{s^{2N - 1} }}{{1 + \left( {st/z} \right)^2 }}ds} } \right)t^{2N} \log \left( {1 - e^{ - 2\pi t} } \right)dt} 
\]
for $\left|\arg z\right|<\frac{\pi}{2}$ and $N\geq 1$. To prove \eqref{eq31}, we perform integration by parts in \eqref{eq5} and employ the integral representation
\[
\mathop{\text{Li}_2} \left(w\right)=  - \int_0^w {\frac{{\log \left( {1 - t} \right)}}{t}dt} ,
\]
which is valid for $\left|\arg\left(w-1\right)\right|<\pi$ \cite[equation 25.12.2]{NIST}.

\subsection{Proof of the formulas \eqref{eq6} and \eqref{eq7}} Our starting point is the functional relationship with the Hurwitz Zeta function $\zeta \left( { s,z} \right)$, due to Vardi \cite{Vardi}:
\begin{gather}\label{eq8}
\begin{split}
\log G\left( {z + 1} \right) & = z\log \Gamma \left( {z + 1} \right) - \zeta '\left( { - 1,z + 1} \right) + \frac{1}{{12}} - \log A
\\ & = z\log \Gamma \left( {z + 1} \right) - \zeta '\left( { - 1,z} \right) - z\log z + \frac{1}{{12}} - \log A,
\end{split}
\end{gather}
for $\Re\left(z\right)>0$, with $\zeta '\left( { s,z} \right) = \partial \zeta \left( {s,z} \right)/\partial s$. Applying formula \cite[equation 25.11.7]{NIST} with $a=z-1$, $n=1$ and the functional equation \cite[equation 25.11.3]{NIST} for the Hurwitz Zeta function, yields
\[
\zeta \left( {s,z} \right) = \frac{{z^{ - s} }}{2} + \frac{{z^{1 - s} }}{{s - 1}} + \frac{{sz^{ - 1 - s} }}{{12}} - \frac{{\Gamma \left( {s + 3} \right)}}{{6\Gamma \left( s \right)}}\int_0^{ + \infty } {\frac{{B_3 \left( {t - \left\lfloor t \right\rfloor } \right)}}{{\left( {t + z} \right)^{s + 3} }}dt} ,
\]
for $z > 1$ and $\Re \left(s \right) >  - 2$ with $s \ne 1$. Differentiation with respect to $s$ and the substitution $s=-1$ gives
\[
\zeta '\left( { - 1,z} \right) =-\frac{1}{4}z^2 + \left( {\frac{z\left(z-1\right)}{2} + \frac{1}{{12}}} \right)\log z +\frac{1}{12} + \frac{1}{6}\int_0^{ + \infty } {\frac{{B_3 \left( {t - \left\lfloor t \right\rfloor } \right)}}{{\left( {t + z} \right)^2 }}dt} .
\]
Substituting this expression into \eqref{eq8}, using the definition of $R_1 \left( z \right)$ and analytic continuation in $z$, we deduce
\begin{equation}\label{eq19}
R_1 \left( z \right) = - \frac{1}{6}\int_0^{ + \infty } {\frac{{B_3 \left( {t - \left\lfloor t \right\rfloor } \right)}}{{\left( {t + z} \right)^2 }}dt} 
\end{equation}
for $\left|\arg z\right| <\pi$. Using the identity
\[
B'_n \left( x \right) = nB_{n - 1} \left( x \right)
\]
and integration by parts, induction on $N$ yields the formulas
\begin{gather}\label{eq9}
\begin{split}
R_N \left( z \right) & =  - \frac{1}{{\left( {2N - 1} \right)2N}}\int_0^{ + \infty } {\frac{{B_{2N} \left( {t - \left\lfloor t \right\rfloor } \right)}}{{\left( {t + z} \right)^{2N - 1} }}dt} \\ & =  - \frac{1}{{2N\left( {2N + 1} \right)}}\int_0^{ + \infty } {\frac{{B_{2N + 1} \left( {t - \left\lfloor t \right\rfloor } \right)}}{{\left( {t + z} \right)^{2N} }}dt} 
\end{split}
\end{gather}
for $N\geq 2$. Thus, the representation \eqref{eq6} is proved. To prove \eqref{eq7} we use the first formula in \eqref{eq9}:
\begin{align*}
R_N \left( z \right)  = \; & \frac{{B_{2N + 2} }}{{2N\left( {2N + 1} \right)\left( {2N + 2} \right)z^{2N} }} + R_{N + 1} \left( z \right)
\\  =  & - \frac{1}{{\left( {2N + 1} \right)\left( {2N + 2} \right)}}\int_0^{ + \infty } {\frac{{B_{2N + 2} }}{{\left( {t + z} \right)^{2N + 1} }}dt} \\ & - \frac{1}{{\left( {2N + 1} \right)\left( {2N + 2} \right)}}\int_0^{ + \infty } {\frac{{B_{2N + 2} \left( {t - \left\lfloor t \right\rfloor } \right)}}{{\left( {t + z} \right)^{2N + 1} }}dt} 
\\ =  & - \frac{1}{{\left( {2N + 1} \right)\left( {2N + 2} \right)}}\int_0^{ + \infty } {\frac{{B_{2N + 2}  + B_{2N + 2} \left( {t - \left\lfloor t \right\rfloor } \right)}}{{\left( {t + z} \right)^{2N + 1} }}dt} .
\end{align*}

\section{Proof of the error bounds}\label{section3} In this section, we prove Theorems \ref{thm2} and \ref{thm3}. To make the subsequent formulas simpler, we introduce the notation
\begin{equation}\label{eq3}
\ell \left( \theta  \right) = \begin{cases} \left|\csc\left(2\theta\right)\right| & \text{ if } \; \frac{\pi}{4} < \left|\theta\right| <\frac{\pi}{2} \\ 1 & \text{ if } \; \left|\theta\right| \leq \frac{\pi}{4} . \end{cases}
\end{equation}
It is easy to show that for any $r>0$,
\begin{equation}\label{eq12}
\frac{1}{\left| {1 + re^{ - 2\theta i} } \right|} \le \ell \left( \theta  \right).
\end{equation}
From now on, we write $\theta = \arg z$. Applying this inequality and the formula \eqref{eq11} for the Bernoulli numbers to the representation \eqref{eq5}, we obtain
\begin{align*}
\left| {R_N \left( z \right)} \right| & \le \frac{1}{{\left| z \right|^{2N} }}\frac{1}{\pi }\int_0^{ + \infty } {\left( {\int_0^1 {\left| {\frac{{s^{2N - 1} }}{{1 + \left( {st/\left| z \right|} \right)^2 e^{ - 2i\theta } }}} \right|ds} } \right)t^{2N} \left| {\log \left( {1 - e^{ - 2\pi t} } \right)} \right|dt} 
\\ & \le \ell \left( \theta  \right)\frac{1}{{\left| z \right|^{2N} }}\frac{1}{\pi }\int_0^{ + \infty } {\left( {\int_0^1 {s^{2N - 1} ds} } \right)t^{2N} \left| {\log \left( {1 - e^{ - 2\pi t} } \right)} \right|dt} 
\\ & = \ell \left( \theta  \right)\frac{1}{{2N\left| z \right|^{2N} }}\frac{{ - 1}}{\pi }\int_0^{ + \infty } {t^{2N} \log \left( {1 - e^{ - 2\pi t} } \right)dt} 
\\ & = \frac{{\left| {B_{2N + 2} } \right|}}{{2N\left( {2N + 1} \right)\left( {2N + 2} \right)\left| z \right|^{2N} }}\ell \left( \theta  \right).
\end{align*}
This partially proves \eqref{eq15}. To prove the statement for the case $z>0$, we note that $0 < \frac{1}{1 + \left( {st/z} \right)^2} < 1$ for any positive $s$, $t$ and $z$, and therefore from \eqref{eq5} we obtain
\begin{align*}
0 < \left( { - 1} \right)^N R_N \left( z \right) & < \frac{1}{{z^{2N} }}\frac{-1}{\pi }\int_0^{ + \infty } {\left( {\int_0^1 {s^{2N - 1} ds} } \right)t^{2N} \log \left( {1 - e^{ - 2\pi t} } \right)dt} \\
& = \left( { - 1} \right)^N \frac{{B_{2N + 2} }}{{2N\left( {2N + 1} \right)\left( {2N + 2} \right)z^{2N} }}.
\end{align*}
Now, let $0 < \varphi  < \frac{\pi }{2}$ be an acute angle that may depend on $N$. First, suppose that $\frac{\pi }{4} + \varphi  < \arg z  < \frac{\pi }{2}+ \varphi$. To obtain the analytic continuation of the formula \eqref{eq5} for $R_N \left( z \right)$ to this sector, we rotate the path of integration in \eqref{eq5} through the angle $\varphi$. Doing so, we find
\begin{align*}
R_N \left( z \right) = \; & \frac{1}{{z^{2N} }}\frac{{\left( { - 1} \right)^{N + 1} }}{\pi }\int_0^{ + \infty e^{i\varphi } } {\left( {\int_0^1 {\frac{{s^{2N - 1} }}{{1 + \left( {st/z} \right)^2 }}ds} } \right)t^{2N} \log \left( {1 - e^{ - 2\pi t} } \right)dt} 
\\ = \; & \frac{1}{{z^{2N} }}\frac{{\left( { - 1} \right)^{N + 1} }}{\pi }\left( {\frac{{e^{i\varphi } }}{{\cos \varphi }}} \right)^{2N + 1} \\ & \times \int_0^{ + \infty } {\left( {\int_0^1 {\frac{{s^{2N - 1} }}{{1 + \left( {st/\left| z \right|\cos \varphi } \right)^2 e^{ - 2i\left( {\theta  - \varphi } \right)} }}ds} } \right)t^{2N} \log \left( {1 - e^{ - 2\pi \frac{{te^{i\varphi } }}{{\cos \varphi }}} } \right)dt} ,
\end{align*}
for $\frac{\pi }{4} + \varphi  < \arg z  < \frac{\pi }{2}+ \varphi$ and $N\geq 1$. To get a simple estimation for $R_N \left( z \right)$, we first note that
\[
\left| {\log \left( {1 - e^{ - 2\pi \frac{{te^{i\varphi } }}{{\cos \varphi }}} } \right)} \right| \leq -\log \left( {1 - \left|e^{ - 2\pi \frac{{te^{i\varphi } }}{{\cos \varphi }}}\right| } \right)  = -\log \left( {1 - e^{ - 2\pi t} } \right)
\]
for any $t>0$ and $0 < \varphi  < \frac{\pi }{2}$. Employing this inequality, the bound \eqref{eq12} and the formula \eqref{eq11} for the Bernoulli numbers, we deduce
\begin{align*}
\left| {R_N \left( z \right)} \right| & \le \frac{{\csc \left( {2\left( {\theta  - \varphi } \right)} \right)}}{{\cos ^{2N + 1} \varphi }}\frac{1}{{\left| z \right|^{2N} }}\frac{-1}{\pi }\int_0^{ + \infty } {\left( {\int_0^1 {s^{2N - 1} ds} } \right)t^{2N} \log \left( {1 - e^{ - 2\pi t} } \right)dt} 
\\ & = \frac{{\csc \left( {2\left( {\theta  - \varphi } \right)} \right)}}{{\cos ^{2N + 1} \varphi }}\frac{1}{{2N\left| z \right|^{2N} }}\frac{{ - 1}}{\pi }\int_0^{ + \infty } {t^{2N} \log \left( {1 - e^{ - 2\pi t} } \right)dt} 
\\ & = \frac{{\csc \left( {2\left( {\theta  - \varphi } \right)} \right)}}{{\cos ^{2N + 1} \varphi }}\frac{{\left| {B_{2N + 2} } \right|}}{{2N\left( {2N + 1} \right)\left( {2N + 2} \right)\left| z \right|^{2N} }}
\end{align*}
for $\frac{\pi }{4} + \varphi  < \arg z  < \frac{\pi }{2}+ \varphi$ and $N\geq 1$. The minimisation of the factor
\begin{equation}\label{eq10}
\frac{ \csc \left( {2\left( {\theta  - \varphi } \right)} \right)}{\cos^{2N + 1 } \varphi },
\end{equation}
can be done by applying a lemma of Meijer's \cite[p. 956]{Meijer}. In our case, Meijer's lemma gives that the minimising value $\varphi = \varphi^\ast$ in \eqref{eq10}, is the unique solution of the implicit equation
\begin{equation}\label{eq14}
\left( {2N + 3 } \right)\cos \left( {3\varphi^\ast - 2\theta} \right) = \left( {2N - 1 } \right)\cos \left( {\varphi^\ast- 2\theta } \right),
\end{equation}
that satisfies $-\frac{\pi}{2}+\theta < \varphi^\ast <\frac{\pi}{2}$ if $\frac{3\pi}{4} \leq \theta <\pi$; $-\frac{\pi}{2}+\theta < \varphi^\ast <-\frac{\pi}{4}+\theta$ if $\frac{\pi}{2} \leq \theta <\frac{3\pi}{4}$; and $0 < \varphi^\ast <-\frac{\pi}{4}+\theta$ if $\frac{\pi}{4} < \theta <\frac{\pi}{2}$. With this choice of $\varphi$,
\begin{equation}\label{eq13}
\left| {R_N \left( z \right)} \right| \le \frac{{\csc \left( {2\left( {\theta  - \varphi^\ast } \right)} \right)}}{{\cos ^{2N + 1} \varphi^\ast }}\frac{{\left| {B_{2N + 2} } \right|}}{{2N\left( {2N + 1} \right)\left( {2N + 2} \right)\left| z \right|^{2N} }}
\end{equation}
for $\frac{\pi}{4}<\arg z <\pi$. By the reflection principle, we have $\left| {R_N \left( {\bar z} \right)} \right| = \left| {\overline {R_N \left( z \right)} } \right| = \left| {R_N \left( z \right)} \right|$, which, in combination with \eqref{eq13}, gives
\[
\left| {R_N \left( z \right)} \right| \le  - \frac{{\csc \left( {2\left( {\theta  - \varphi ^ *  } \right)} \right)}}{{\cos ^{2N + 1} \varphi ^ *  }}\frac{{\left| {B_{2N + 2} } \right|}}{{2N\left( {2N + 1} \right)\left( {2N + 2} \right)\left| z \right|^{2N} }}
\]
for $- \pi  < \arg z <  - \frac{\pi }{4}$. Here $-\frac{\pi }{2} < \varphi^\ast  < 0$ is the unique solution of the equation \eqref{eq14}, that satisfies $ -\frac{\pi}{2} < \varphi^\ast < \frac{\pi}{2}+\theta$ if $-\pi <\theta \leq -\frac{3\pi}{4}$; $\frac{\pi}{4}+\theta < \varphi^\ast < \frac{\pi}{2}+\theta$ if $ -\frac{3\pi}{4} < \theta  \leq -\frac{\pi}{2}$; and $\frac{\pi}{4}+\theta < \varphi^\ast < 0$ if $-\frac{\pi}{2} < \theta <-\frac{\pi}{4}$. The proof of Theorem \ref{thm3} is complete.

To complete the proof of \eqref{eq15}, note that when $\arg z = \theta = \frac{\pi}{2}$, the minimising value $\varphi^\ast$ in \eqref{eq14} is given explicitly by
\[
\varphi^\ast  = \arctan \left( {\frac{1}{{\sqrt {2N + 2} }}} \right),
\]
and therefore we have
\begin{align*}
\frac{{\csc \left( {2\left( {\theta  - \varphi^\ast } \right)} \right)}}{{\cos^{2N + 1} \varphi^\ast }} \le \frac{{\csc \left( {2\left( {\frac{\pi }{2} - \varphi^\ast } \right)} \right)}}{{\cos^{2N + 1 } \varphi^\ast }} & = \frac{1}{2}\left( {1 + \frac{1}{{2N + 2 }}} \right)^{N} \sqrt {2N + 2} \\ & \le \frac{1}{2}\sqrt {e\left( {2N + \frac{5}{2}} \right)}, 
\end{align*}
as long as $\frac{\pi }{4} + \varphi^\ast < \theta  \leq \frac{\pi }{2}$ and $N\geq 1$. We also have
\[
\frac{1}{2}\sqrt {e\left( {2N + \frac{5}{2}} \right)}  \ge \frac{1}{2}\sqrt {e\left( {2 + \frac{5}{2}} \right)}  \ge \csc \left( {2\theta } \right),
\]
for $\frac{\pi }{4} < \theta  \le \frac{\pi }{4} + \varphi ^\ast \le \frac{\pi }{4} + \arctan \left( {\frac{1}{2}} \right)$, therefore
\[
\left| {R_N \left( z \right)} \right| \le \frac{1}{2}\sqrt {e\left( {2N + \frac{5}{2}} \right)}\frac{{\left| {B_{2N + 2} } \right|}}{{2N\left( {2N + 1} \right)\left( {2N + 2} \right)\left| z \right|^{2N} }},
\]
provided that $\frac{\pi }{4} < \theta  \leq \frac{\pi }{2}$ and $N\geq 1$. For the conjugate sector $-\frac{\pi }{2} < \theta  \leq -\frac{\pi }{4}$, we can simply use $\left| {R_N \left( {\bar z} \right)} \right| = \left| {\overline {R_N \left( z \right)} } \right| = \left| {R_N \left( z \right)} \right|$. This finishes the proof of \eqref{eq15}.

Finally, we prove \eqref{eq17} with the help of the expression \eqref{eq7}. First, we show that the quantity $B_{2N + 2}  + B_{2N + 2} \left( {t - \left\lfloor t \right\rfloor } \right)$ does not change sign for any $t\geq 0$ with a fixed $N\geq 1$. Indeed, by the Fourier series of the Bernoulli polynomials \cite[equation 24.8.1]{NIST} we have
\[
\left( { - 1} \right)^N \left( {B_{2N + 2}  + B_{2N + 2} \left( {t - \left\lfloor t \right\rfloor } \right)} \right) = \frac{{2\left( {2N + 2} \right)!}}{{\left( {2\pi } \right)^{2N + 2} }}\sum\limits_{k = 1}^\infty  {\frac{{1 + \cos \left( {2\pi kt} \right)}}{{k^{2N + 2} }}}  \ge 0.
\]
It is easy to show that
\[
\left| {t + z} \right| \ge \left( {t + \left| z \right|} \right)\cos \left( {\frac{\theta }{2}} \right)
\]
for $\left|\arg z\right|<\pi$ and $t\geq 0$. Applying this inequality for \eqref{eq7}, we obtain the desired result
\begin{align*}
\left| {R_N \left( z \right)} \right| & \le \frac{1}{{\left( {2N + 1} \right)\left( {2N + 2} \right)}}\int_0^{ + \infty } {\frac{{\left| {B_{2N + 2}  + B_{2N + 2} \left( {t - \left\lfloor t \right\rfloor } \right)} \right|}}{{\left| {t + z} \right|^{2N + 1} }}dt} 
\\ & = \left| {\frac{1}{{\left( {2N + 1} \right)\left( {2N + 2} \right)}}\int_0^{ + \infty } {\frac{{B_{2N + 2}  + B_{2N + 2} \left( {t - \left\lfloor t \right\rfloor } \right)}}{{\left| {t + z} \right|^{2N + 1} }}dt} } \right|
\\ & \le \left| {\frac{1}{{\left( {2N + 1} \right)\left( {2N + 2} \right)}}\int_0^{ + \infty } {\frac{{B_{2N + 2}  + B_{2N + 2} \left( {t - \left\lfloor t \right\rfloor } \right)}}{{\left( {t + \left| z \right|} \right)^{2N + 1} }}dt} } \right|\sec ^{2N + 1} \left( {\frac{\theta }{2}} \right)
\\ & = \left| {R_N \left( {\left| z \right|} \right)} \right|\sec ^{2N + 1} \left( {\frac{\theta }{2}} \right) \le \frac{{\left| {B_{2N + 2} } \right|}}{{2N\left( {2N + 1} \right)\left( {2N + 2} \right)\left| z \right|^{2N} }}\sec ^{2N + 1} \left( {\frac{\theta }{2}} \right).
\end{align*}
In the final step we have made use of \eqref{eq15}.

\section{Exponentially improved asymptotic expansion}\label{section4}

We shall find it convenient to express our exponentially improved expansion in terms of the (scaled) Terminant function, which is defined in terms of the Incomplete gamma function as
\[
\widehat T_p \left( w \right) = \frac{{e^{\pi ip} \Gamma \left( p \right)}}{{2\pi i}}\Gamma \left( {1 - p,w} \right) = \frac{e^{\pi ip} w^{1 - p} e^{ - w} }{2\pi i}\int_0^{ + \infty } {\frac{{t^{p - 1} e^{ - t} }}{w + t}dt} \; \text{ for } \; p>0 \; \text{ and } \; \left| \arg w \right| < \pi ,
\]
and by analytic continuation elsewhere. Olver \cite{Olver1} showed that when $p \sim \left|w\right|$ and $w \to \infty$, we have
\begin{equation}\label{eq27}
\widehat T_p \left( w \right) = \frac{1}{2} + \frac{1}{2}\mathop{\text{erf}} \left( {c\left( \varphi  \right)\sqrt {\frac{1}{2}\left| w \right|} } \right) + \mathcal{O}\left( {\frac{{e^{ - \frac{1}{2}\left| w \right|c^2 \left( \varphi  \right)} }}{{\left| w \right|^{\frac{1}{2}} }}} \right)
\end{equation}
for $-\pi +\delta \leq \arg w \leq 3 \pi -\delta$, $0 < \delta  \le 2\pi$; and
\begin{equation}\label{eq28}
e^{ - 2\pi ip} \widehat T_p \left( w \right) =  - \frac{1}{2} + \frac{1}{2}\mathop{\text{erf}} \left( { - \overline {c\left( { - \varphi } \right)} \sqrt {\frac{1}{2}\left| w \right|} } \right) + \mathcal{O}\left( {\frac{{e^{ - \frac{1}{2}\left| w \right|\overline {c^2 \left( { - \varphi } \right)} } }}{{\left| w \right|^{\frac{1}{2}} }}} \right)
\end{equation}
for $- 3\pi  + \delta  \le \arg w \le \pi  - \delta$, $0 < \delta \le 2\pi$. Here $\varphi = \arg w$ and $\mathop{\text{erf}}$ denotes the Error function. The quantity $c\left( \varphi  \right)$ is defined implicitly by the equation
\[
\frac{1}{2}c^2 \left( \varphi  \right) = 1 + i\left( {\varphi  - \pi } \right) - e^{i\left( {\varphi  - \pi } \right)},
\]
and corresponds to the branch of $c\left( \varphi  \right)$ which has the following expansion in the neighbourhood of $\varphi = \pi$:
\begin{equation}\label{eq29}
c\left( \varphi  \right) = \left( {\varphi  - \pi } \right) + \frac{i}{6}\left( {\varphi  - \pi } \right)^2  - \frac{1}{{36}}\left( {\varphi  - \pi } \right)^3  - \frac{i}{{270}}\left( {\varphi  - \pi } \right)^4  +  \cdots .
\end{equation}
For complete asymptotic expansions, see Olver \cite{Olver1}. We remark that Olver uses the different notation $F_p \left( w \right) = ie^{ - \pi ip} \widehat T_p \left( w \right)$ for the Terminant function and the other branch of the function $c\left( \varphi  \right)$. For further properties of the Terminant function, see, for example, Paris and Kaminski \cite[Chapter 6]{Paris}.

\subsection{The discontinuous treatment of the Stokes multipliers}\label{subsection41} The asymptotic expansion of $\log G\left(z+1\right)$ in the sector $\left|\arg z\right|<\frac{\pi}{2}$ is given by \eqref{eq4}. We shall obtain an expansion in the sector $\frac{\pi}{2}<\left|\arg z\right|<\pi$. We expand the Dilogarithm in \eqref{eq31} using its series \eqref{eq18} and apply partial fraction decomposition to deduce
\begin{equation}\label{eq30}
R_1 \left( z \right)  =  \sum\limits_{k = 1}^\infty  {\frac{1}{{2\pi k^2 }}\left( {\frac{1}{{2\pi i}}\int_0^{ + \infty } {\frac{{e^{ - 2\pi kt} }}{{it - z}}dt}  + \frac{1}{{2\pi i}}\int_0^{ + \infty } {\frac{{e^{ - 2\pi kt} }}{{it + z}}dt} } \right)}. 
\end{equation}
To find the analytic continuation of $R_1 \left( z \right)$ to the sector $\frac{\pi}{2}<\arg z < \pi$, we rotate the path of integration in the first integral to obtain
\[
R_1 \left( z \right) =  -\sum\limits_{k = 1}^\infty  {\frac{{e^{2\pi ikz} }}{{2\pi i k^2 }}}  + R_1 \left( {ze^{ - \pi i} } \right),
\]
for $\frac{\pi}{2}<\arg z < \pi$. By the connection formula $R_1 \left( {\bar z} \right)= \overline {R_1 \left( z \right)}$, we have, similarly,
\[
R_1 \left( z \right) = \sum\limits_{k = 1}^\infty  {\frac{{e^{ - 2\pi ikz} }}{{2\pi i k^2 }}}  + R_1 \left( {ze^{\pi i} } \right),
\]
for $-\pi<\arg z <-\frac{\pi}{2}$. Using the asymptotic expansion \eqref{eq4} for $R_1 \left( {ze^{ \mp \pi i} } \right)$, we finally obtain
\begin{equation}\label{eq21}
R_1 \left( z \right) \sim  \mp \sum\limits_{k = 1}^\infty  {\frac{{e^{ \pm 2\pi ikz} }}{{2\pi i k^2 }}}  + \sum\limits_{n = 1}^\infty  {\frac{{B_{2n + 2} }}{{2n\left( {2n + 1} \right)\left( {2n + 2} \right)z^{2n} }}} 
\end{equation}
as $z\to \infty$ in the sectors $\frac{\pi}{2}< \pm \arg z < \pi$.

The case $\left|\arg z\right| = \frac{\pi}{2}$ must be considered separately. On the positive imaginary axis we put $z=iy$, $y>0$. Integrating once by parts in \eqref{eq19} and using the Fourier series of the Bernoulli polynomial $B_2 \left( {t - \left\lfloor t \right\rfloor } \right)$ \cite[equation 24.8.1]{NIST}, we find
\begin{align*}
R_1 \left( {iy} \right) & =  - \frac{1}{2}\int_0^{ + \infty } {\frac{{B_2 \left( {t - \left\lfloor t \right\rfloor } \right)}}{{t + iy}}dt}  = \frac{{iy}}{2}\int_0^{ + \infty } {\frac{{B_2 \left( {t - \left\lfloor t \right\rfloor } \right)}}{{t^2  + y^2 }}dt}  - \frac{1}{2}\int_0^{ + \infty } {\frac{{tB_2 \left( {t - \left\lfloor t \right\rfloor } \right)}}{{t^2  + y^2 }}dt} 
\\ & = \frac{{2iy}}{{\left( {2\pi } \right)^2 }}\sum\limits_{n = 1}^\infty  {\frac{1}{{n^2 }}\int_0^{ + \infty } {\frac{{\cos \left( {2\pi nt} \right)}}{{t^2  + y^2 }}dt} }  - \frac{2}{{\left( {2\pi } \right)^2 }}\sum\limits_{n = 1}^\infty  {\frac{1}{{n^2 }}\int_0^{ + \infty } {\frac{{t\cos \left( {2\pi nt} \right)}}{{t^2  + y^2 }}dt} } .
\end{align*}
The first integral can be evaluated using Fourier transforms, and we find
\[
\frac{{2iy}}{{\left( {2\pi } \right)^2 }}\sum\limits_{n = 1}^\infty  {\frac{1}{{n^2 }}\int_0^{ + \infty } {\frac{{\cos \left( {2\pi nt} \right)}}{{t^2  + y^2 }}dt} }  =  - \sum\limits_{n = 1}^\infty  {\frac{1}{2} \frac{{e^{ - 2\pi ny} }}{{2\pi i n^2 }}} .
\]
The second integral can be written in terms of the Exponential integral \cite[equation 3.723/5.]{Tables}, hence we arrive at
\begin{equation}\label{eq20}
 - \frac{2}{{\left( {2\pi } \right)^2 }}\sum\limits_{n = 1}^\infty  {\frac{1}{{n^2 }}\int_0^{ + \infty } {\frac{{t\cos \left( {2\pi nt} \right)}}{{t^2  + y^2 }}dt} }  = \frac{1}{{\left( {2\pi } \right)^2 }}\sum\limits_{\substack{k =  - \infty \\ k \ne 0}}^\infty  {\frac{1}{{k^2 }}e^{ - 2\pi ky} \mathop{\text{Ei}} \left( {2\pi ky} \right)} .
\end{equation}
From Watson's lemma \cite[pp. 71--72]{Olver3}\cite[pp. 20--21]{Wong}, we obtain that
\[
e^{ \mp 2\pi ky} \mathop{\text{Ei}} \left( { \pm 2\pi ky} \right) =  \pm \mathrm{PV} \int_0^{ + \infty } {\frac{{e^{ - 2\pi kyt} }}{{1 \mp t}}dt \sim \sum\limits_{n = 1}^\infty  {\frac{{\left( {n - 1} \right)!}}{{\left( { \mp 2\pi ky} \right)^n }}} } ,\; k > 0,
\]
as $y \to +\infty$, where $\mathrm{PV}$ denotes the Cauchy principal value. Substituting this expansion into \eqref{eq20}, changing the order of summation and using the representation of the Bernoulli numbers in terms of the Riemann Zeta function \cite[equation 25.6.2]{NIST}, we find
\[
 - \frac{2}{{\left( {2\pi } \right)^2 }}\sum\limits_{n = 1}^\infty  {\frac{1}{{n^2 }}\int_0^{ + \infty } {\frac{{t\cos \left( {2\pi nt} \right)}}{{t^2  + y^2 }}dt} }  \sim \sum\limits_{n = 1}^\infty  {\frac{{B_{2n + 2} }}{{2n\left( {2n + 1} \right)\left( {2n + 2} \right)\left( {iy} \right)^{2n} }}} ,
\]
as $y \to +\infty$. Combining all the partial results, we arrive at
\begin{equation}\label{eq22}
R_1 \left( z \right) \sim  - \sum\limits_{k = 1}^\infty  {\frac{1}{2} \frac{{e^{2\pi ikz} }}{{2\pi i k^2 }}}  + \sum\limits_{n = 1}^\infty  {\frac{{B_{2n + 2} }}{{2n\left( {2n + 1} \right)\left( {2n + 2} \right)z^{2n} }}} 
\end{equation}
as $z\to \infty$ with $\arg z = \frac{\pi}{2}$. By the connection formula $R_1 \left( {\bar z} \right)= \overline {R_1 \left( z \right)}$, we have, similarly,
\begin{equation}\label{eq23}
R_1 \left( z \right) \sim  \sum\limits_{k = 1}^\infty  {\frac{1}{2} \frac{{e^{-2\pi ikz} }}{{2\pi i k^2 }}}  + \sum\limits_{n = 1}^\infty  {\frac{{B_{2n + 2} }}{{2n\left( {2n + 1} \right)\left( {2n + 2} \right)z^{2n} }}} 
\end{equation}
as $z\to \infty$ with $\arg z = -\frac{\pi}{2}$. We can unify the series \eqref{eq4}, \eqref{eq21}, \eqref{eq22} and \eqref{eq23} as follows:
\begin{gather}\label{eq24}
\begin{split}
\log G\left( {z + 1} \right) \sim \frac{1}{4}z^2  + z\log \Gamma \left( {z + 1} \right) - \left( {\frac{1}{2}z\left( {z + 1} \right) + \frac{1}{{12}}} \right)\log z - \log A \\ +\sum\limits_{n = 1}^\infty  {S_k \left( \theta  \right)e^{ \pm 2\pi ikz} }  + \sum\limits_{n = 1}^\infty  {\frac{{B_{2n + 2} }}{{2n\left( {2n + 1} \right)\left( {2n + 2} \right)z^{2n} }}} ,
\end{split}
\end{gather}
where
\begin{equation}\label{eq25}
S_k\left( \theta  \right) = \begin{cases} 0 & \text{ if }  \left|\theta\right| < \frac{\pi}{2} \\ \mp \cfrac{1}{2}\cfrac{1}{{2\pi ik^2 }} & \text{ if } \theta = \pm \frac{\pi}{2} \\ \mp \cfrac{1}{{2\pi ik^2 }} & \text{ if } \frac{\pi}{2} < \left|\theta\right| < \pi ,\end{cases}
\end{equation}
and $\theta = \arg z$. The upper or lower sign is taken in \eqref{eq24} and \eqref{eq25} according as $z$ is in the upper or lower half-plane. It seems that there is a discontinuous change in the coefficients of the exponential terms when $\arg z$ changes continuously across $\arg z = \pm \frac{\pi}{2}$. We have encountered a Stokes phenomenon with Stokes lines $\arg z = \pm \frac{\pi}{2}$. The function $S_k \left( \theta  \right)$ is called the Stokes multiplier of the subdominant exponential $e^{ \pm 2\pi ikz}$. The forms of these multipliers are in agreement with Dingle's non-rigorous ``final main rule" in his theory of terminants \cite[p. 414]{Dingle}, namely that half the discontinuity occurs on reaching the Stokes line, and half on leaving it the other side.

In the important papers \cite{Berry1}\cite{Berry2}, Berry provided a new interpretation of the Stokes phenomenon; he found that assuming optimal truncation, the transition between compound asymptotic expansions is of Error function type, thus yielding a smooth, although very rapid, transition as a Stokes line is crossed.

Using the expansion given in Theorem \ref{thm4}, we shall show in the next subsection that the asymptotic expansion of the Barnes $G$-function exhibits the Berry transition between the two asymptotic series across the Stokes lines $\arg z = \pm \frac{\pi}{2}$. More precisely, we shall find that the subdominant exponentials $e^{ \pm 2\pi i k z}$ ``switch on" in a rapid and smooth way as $\arg z$ passes through the values $\pm\frac{\pi}{2}$.

\subsection{Smoothing the Stokes discontinuities}\label{subsection42} We begin with the proof of Theorem \ref{thm4}. We employ the well-known expression for non-negative integer $N_k$
\[
\frac{1}{1 - x} = \sum\limits_{n = 0}^{2N_k-1} {x^n}  + \frac{x^{2N_k}}{1 - x},\; x \neq 1,
\]
to expand the fractions under the integrals in \eqref{eq30} in powers of $\pm it/z$, leading to the result
\begin{align*}
R_1 \left( z \right) =  & -\sum\limits_{k = 1}^\infty  {\frac{1}{{\left( {2\pi k} \right)^2 }}\sum\limits_{n = 0}^{N_k-1} {\left( { - 1} \right)^n \frac{{2\left( {2n + 1} \right)!}}{{\left( {2\pi kz} \right)^{2n+2} }}} } 
\\ & + \sum\limits_{k = 1}^\infty  {\frac{1}{{2\pi k^2 }}\left( {-\frac{{\left( { - 1} \right)^{N_k } }}{{z^{2N_k } }}\frac{1}{{2\pi i}}\int_0^{ + \infty } {\frac{{t^{2N_k } e^{ - 2\pi kt} }}{{z - it}}dt}  + \frac{{\left( { - 1} \right)^{N_k } }}{{z^{2N_k } }}\frac{1}{{2\pi i}}\int_0^{ + \infty } {\frac{{t^{2N_k } e^{ - 2\pi kt} }}{{z + it}}dt} } \right)} ,
\end{align*}
for $\left|\arg z\right|<\frac{\pi}{2}$. The integrals can be written in terms of the Terminant function, and by analytic continuation, we arrive at
\begin{align*}
R_1 \left( z \right) = & - \sum\limits_{k = 1}^\infty  {\frac{1}{{\left( {2\pi k} \right)^2 }}\sum\limits_{n = 0}^{N_k  - 1} {\left( { - 1} \right)^n \frac{{2\left( {2n + 1} \right)!}}{{\left( {2\pi kz} \right)^{2n + 2} }}} } \\ & - \sum\limits_{k = 1}^\infty  {\left( {\frac{{\widehat T_{2N_k  + 1} \left( {2\pi kiz} \right)}}{{2\pi ik^2 }}e^{2\pi kiz}  + \frac{{\widehat T_{2N_k  + 1} \left( { - 2\pi kiz} \right)}}{{2\pi ik^2 }}e^{ - 2\pi kiz} } \right)} ,
\end{align*}
for $\left|\arg z\right|<\pi$. The absolute convergence of the infinite series involving the Terminant functions is justified by the large $w$ behaviour of $\widehat T_p \left( w \right)$ (see, e.g., Paris and Kaminski \cite[p. 260]{Paris}). The proof of Theorem \ref{thm4} is complete.

Now, we show that the asymptotic expansion of the Barnes $G$-function has the Berry transition property. We choose $N_k \sim \pi k \left| z\right|$ for each positive $k$. With this choice, in the upper half-plane the terms involving $\widehat T_{2N_k  + 1} \left( { - 2\pi kiz} \right)$ are exponentially small, the dominant contribution comes from the terms involving $\widehat T_{2N_k  + 1} \left( {2\pi kiz} \right)$. From \eqref{eq27} and \eqref{eq29}, the Terminant functions have the asymptotic behaviour
\[
\widehat T_{2N_k  + 1} \left( {2\pi kiz} \right) \sim \frac{1}{2} + \frac{1}{2}\mathop{\text{erf}} \left( {\left( {\theta  - \frac{\pi }{2}} \right)\sqrt {\pi k\left| z \right|} } \right)
\]
provided that $\arg z = \theta$ is close to $\frac{\pi}{2}$ and $z$ is large. Similarly, in the lower half-plane, the dominant contribution is
controlled by the terms involving $\widehat T_{2N_k  + 1} \left( { - 2\pi kiz} \right)$. From \eqref{eq28} and \eqref{eq29}, we have
\[
\widehat T_{2N_k  + 1} \left( { - 2\pi kiz} \right) \sim - \frac{1}{2} + \frac{1}{2}\mathop{\text{erf}} \left( {\left( {\theta  + \frac{\pi }{2}} \right)\sqrt {\pi k\left| z \right|} } \right)
\]
under the assumptions that $\arg z = \theta$ is close to $-\frac{\pi}{2}$ and $z$ is large. Therefore, the approximate functional form of the Stokes multipliers near $\arg z = \pm\frac{\pi}{2}$ is
\[
\mp \frac{1}{{2\pi ik^2 }}\left( {\frac{1}{2} \pm \frac{1}{2}\mathop{\text{erf}} \left( {\left( {\theta  \mp \frac{\pi }{2}} \right)\sqrt {\pi k\left| z \right|} } \right)} \right).
\]
For example, in the upper half-plane, the Stokes multipliers change from approximately $0$ when $\theta < \frac{\pi}{2}$ to approximately $-\frac{1}{2\pi ik^2}$ when $\theta > \frac{\pi}{2}$. On the Stokes line $\theta = \frac{\pi}{2}$ they take the value $-\frac{1}{2}\frac{1}{2\pi ik^2}$ up to an exponentially small error. Hence, we have proved that the transition of the Stokes multipliers across the Stokes lines $\arg z = \pm\frac{\pi}{2}$ is smooth.

\section{Discussion}\label{section5} In this paper we have established new integral representations for the remainder term of the asymptotic expansion of the logarithm of the Barnes $G$-function. Using these representations, we have derived new simple, explicit and numerically computable error bounds for this asymptotic series. We have found that along the imaginary axis, suddenly infinitely many exponentially small terms appear in the asymptotic expansion of the Barnes $G$-function. By employing one of the representations for the remainder term, we have obtained an exponentially improved asymptotic expansion for the logarithm of the Barnes $G$-function, which shows that the appearance of these exponentially small terms is in fact smooth, showing the Berry transition property of the asymptotic series of the $G$-function.

If we substitute into \eqref{eq4} the well-known asymptotic expansion of $\log \Gamma \left( {z + 1} \right)$, we obtain
\begin{equation}\label{eq33}
\log G \left( {z + 1} \right) \sim  - \frac{3}{4}z^2  + \frac{z}{2}\log \left( {2\pi } \right) + \left( {\frac{1}{2}z^2  - \frac{1}{{12}}} \right)\log z + \frac{1}{{12}} - \log A + \sum\limits_{n = 1}^\infty  {\frac{{B_{2n + 2} }}{{2n\left( {2n + 2} \right)z^{2n} }}} ,
\end{equation}
as $z\to \infty$ in the sector $\left| {\arg z} \right| \leq \pi -\delta< \pi$ with any fixed $0<\delta\leq \pi$. This is exactly the asymptotic series that was given originally by Barnes \cite{Barnes2}. Using the known error bounds \cite{FWS} and the exponentially improved version \cite[pp. 279--288]{Paris} of the asymptotic series of $\log \Gamma \left( {z + 1} \right)$, together with the results of this paper, similar properties of Barnes' series \eqref{eq33} can be obtained.

\end{document}